\documentclass{amsart}
\usepackage{amsmath}%
\usepackage{amsfonts}%
\usepackage{amssymb}%
\usepackage{graphicx}
\newtheorem{theorem}{Theorem}
\theoremstyle{plain}

\newtheorem{corollary}{Corollary}

\newtheorem{lemma}{Lemma}

\numberwithin{equation}{section}
\begin{document}
\title[On the distribution of coefficients of polynomials]{On the distribution of coefficients of residue polynomials}
\author{L\'aszl\'o Major}
\address{L\'aszl\'o Major\newline%
\indent Institute of Mathematics,   \newline%
\indent Tampere University of Technology,  \newline%
\indent PL 553, 33101 Tampere, Finland }\email{laszlo.major@tut.fi}
\date{July 1, 2010}

\begin{abstract}
Using the formalism of polynomials with positive coefficients, the fact that exactly half of all subsets of a finite set have even cardinality can be generalized asymptotically.
\end{abstract}
\maketitle
The well-known fact that every finite nonempty set has as many subsets of even cardinality as of odd cardinality can be restated as follows. If $d$ denotes the particular integer $2$, then for any integers $0\leq j\leq d-1$, $n\geq 1$ and the standard expansion $c_{n,0}+c_{n,1}x+c_{n,2}x^2+\cdots$ of the polynomial $(\frac 12 + \frac 12 x)^n$ we have
$$\sum_{k\equiv j \hspace{-2mm}\mod \hspace{-0mm}d}\hspace{-2mm}c_{n,k}=\frac 1d.$$
For $n=0$ this obviously fails. It also fails in general if $d\geq 3.$ We will show however, that it remains asymptotically true for all $d\geq 2$ and not only for the polynomial $\frac 12 + \frac 12 x$, but for any polynomial with positive coefficients summing to $1$ (Theorem \ref{egy}). In particular, this will imply that summing up every $d^{th}$ entry in the $n^{th}$ row in Pascal's triangle asymptotically yields $2^n/d.$
 
Questions of distribution about powers of polynomials with positive coefficients, in connection with binomial and multinomial coefficients in particular, have been studied both by classical and contemporary authors. When Euler (1765) was investigating the properties of the trinomial coefficients (see Andrews \cite{andras}), he obtained an unimodal distribution of the coefficients by expanding the $n^{th}$ power of the polynomial $1+x+x^{2}$. In the analogous case of the polynomial $1+x$ the corresponding unimodal distribution is the $n^{th}$ row of Pascal's triangle. In general the distribution of the coefficients of the $n^{th}$ power of the polynomial $p=p_0+p_1x+\cdots+p_mx^{m}$ with nonnegative real coefficients is not necessarily unimodal, but a sufficient condition is given  for unimodality by Boros and Moll \cite{boros}. In this paper we are concerned with a question of equidistribution of coefficients. The main result is Theorem \ref{e}. The lemmas and corollary provide background and an extension.

Let $d \geq 2$ be  a fixed positive integer. Let $r(p)$ denote the remainder of the division of the  polynomial $p$ by $x^d-1$, i.e. that unique polynomial of degree less than $d$ for which the polynomial $p$ is congruent to $r(p)$ modulo $x^d-1$. A polynomial with positive (nonnegative) real coefficients is said to be \emph{positive (nonnegative)} polynomial. The following Lemma provides a 
sufficient condition for the powers of a nonnegative polinomial $p$ to have positive remainder by $x^d-1.$
\begin{lemma}\label{egy}
Let $p\in \mathbb{R}[x]$ be a nonnegative, nonzero polynomial. For any $q\in \mathbb{R}[x]$ and $i\in  \mathbb{N}$ let $q_i$ denote the coefficient of $x^i$ in $q$. If  there exist $k,l\in \{0,\ldots,d-1\}$ such that $r(p)_k>0$, $r(p)_l>0$ and gcd$(d,k-l)=1$, then 
 $r(p^{n})$ is a positive polynomial for all integer $n>d-1$. 
\end{lemma}

\begin{proof}
It can be assumed that $k<l$. Let us denote the difference $l-k$ by $h$. It is easy to  see that if $r((x^k+x^l)^{d-1})$ is positive then also $r(p^{d-1})$ is positive. In addition,  $r((x^k+x^l)^{d-1})=r(x^{k(d-1)}(1+x^h)^{d-1})$ is positive if $r((1+x^h)^{d-1})$ is positive, so let us concentrate on the polynomial $r((1+x^h)^{d-1})$. We expand the expression $(1+x^h)^{d-1}$ in the polynomial ring $\mathbb{R}[x]$ by using the Binomial Theorem: 
\begin{equation}\label{eq10}
(1+x^{h})^{d-1}=\sum_{i=0}^{d-1}\tbinom {d-1}{i}(x^{hi}).
\end{equation}
Let us assume that the following congruence holds for some $i,j\in \{0,\ldots,d-1\}$: 
\begin{equation}\label{eq20}
\tbinom {d-1}{i}(x^{hi})\equiv \tbinom {d-1}{j}(x^{hj}) \mod x^d-1. 
\end{equation}
From the equation (\ref{eq20}) follows that ${hi}-{hj}$ is divisible by $d$. But it was assumed that $gcd(h,d)=1$, consequently $i=j$. Therefore $r((1+x^h)^{d-1})$ has exactly $d$ different nonzero coefficients,  in other words it is positive and consequently $r(p^{d-1})$ is also positive.

It can be shown without difficulty that for any $n\geq d-1$ the remainder $r(p^n)$ is positive if $r(p^{d-1})$ is positive.
\end{proof}
In order to be able to deal more efficiently with powers of  polynomials we need the concept and some useful properties of \emph{circulant matrices}. For a classical reference see e.g. Davis \cite{davis}. 
  Let $v=(v_0,v_1,\ldots,v_{d-1})$  be a row vector in $\mathbb{R}^{d}$. The permutation $\rho:\mathbb{R}^d\rightarrow \mathbb{R}^d$ given by
  $$  \rho (v_0,v_1,\ldots,v_{d-1})=(v_{d-1},v_0,\ldots,v_{d-2})$$ is called \emph{cyclic permutation}. The \emph{circulant matrix} associated to the vector $v$ is the $d\times d$ matrix whose $i^{th}$ row is $\rho^{i-1}(v)$, $i=1,\ldots,d$ and it is denoted by
  $$C=circ(v_0,v_1,\ldots,v_{d-1})=circ(v).$$
The product of two circulant matrices is circulant, therefore any positive integer power of a circulant matrix is circulant. We shall use the connection between powers of circulant matrices and the powers of residue polynomials. If $a=(a_0,a_1,\ldots,a_{d-1})$ is the coefficient vector of some polynomial $r(p)\in \mathbb{R}[x]$, then 
   the first row of $(circ(a))^n$ is the coefficient vector of $r(p^n)$. Lally and Fitzpatrick \cite{alstu} give  a general overview in this subject.
  
 A matrix is called \emph{positive (nonnegative)} if all  its entries are  positive (nonnegative) real numbers.
A $d \times d$ nonnegative matrix $M$ is said to be \emph{doubly stochastic} if the sum of the entries in each row and in each column equals $1$. The product of doubly stochastic matrices is doubly stochastic. The following lemma follows from a general result on eigenvectors of positive matrices,  Theorem 8.2.8 in Horn and Johnson \cite{horn}. Here we provide a direct proof.
 A matrix whose entries are all ones will be denoted by $U$.

\begin{lemma}\label{lem}If $M$ is a $d \times d$ positive doubly stochastic matrix, then
\begin{equation}\label{eq30}
\lim_{n\rightarrow \infty}M^n=\frac 1dU.
\end{equation}
\end{lemma}
\begin{proof}
Let us denote the doubly stochastic  matrix $\frac 1dU$ by $J$. If $M=J$, then the statement is trivial. Otherwise the matrix $M$ can be given in the following form:
\begin{equation}\label{eq40}
M=\lambda J+(1-\lambda )M_0,
\end{equation}
where $0<\lambda <1$ and $M_0$ is a doubly stochastic matrix. Let us write the $n^{th}$ power of $M$ by using the Binomial Theorem and the fact that $JD=DJ=J$ for all doubly stochastic matrix $D$,
\[\begin{split}
M^n&=\sum_{i=0}^{n}\binom ni (\lambda J)^{n-i}((1-\lambda )M_0)^i
=(1-\lambda)^nM_0^n+J\sum_{i=0}^{n-1}\binom ni \lambda^{n-i}(1-\lambda )^i\\&=(1-\lambda)^nM_0^n+J(1-(1-\lambda )^{n})=(1-\lambda)^n(M_0^n-J)+J\\
\end{split}\]
From the above expression follows our statement, because $0<1-\lambda <1$ and consequently $lim_{n\rightarrow \infty} (1-\lambda )^{n}=0$.
\end{proof}
In order that the Lemma \ref{lem} could  be applied for the powers of polynomials, we shall restrict our attention to  polynomials $p$ whose coefficients sum to 1 $(p(1)=1)$.  We may do this without loss of generality, because all polynomials can be written in the form  $p=p(1)\cdot p'$, where $p'$ has the mentioned property. It can also be said that the coefficients of $p'$ form a \emph{stochastic vector}. Clearly the sum of the coefficients of $(p')^n$ is also 1 and the sum of the coefficients of $p^n$ is equal to $p(1)^n$. By using Lemma \ref{egy}, Lemma \ref{lem} and the connection between powers of circulant matrices and the powers of residue polynomials we obtain the following theorem:
\begin{theorem}\label{e}
Let $p$ be a nonconstant polynomial with positive real coefficients summing to 1 and $d\geq 2$ a fixed integer. Then for all $0\leq j\leq d-1$ 
\begin{equation}\label{eq50}
 \lim_{n\rightarrow \infty}\hspace{-0mm}\sum_{k\equiv j \hspace{-2mm}\mod \hspace{-0mm}d}\hspace{-2mm}(p^n)_k=\frac 1d
\end{equation}
where $(p^n)_0+ (p^n)_1x+(p^n)_2x^2+\cdots$ is the standard expansion of $p^n$. \hspace{21mm}$\square$
 \end{theorem}
 
 The restriction that $p$ is positive is quite strong.  Theorem \ref{e} can be extended to larger classes of polynomials. For instance the following corollary applies Lemma \ref{egy} so as to provide such an extension.
 
\begin{corollary}\label{c1}Let $d\geq 2$ a fixed integer. If the nonnegative polynomial $p$ satisfies the conditions of Lemma \ref{egy} and its coefficients sum to 1, then the convergence (\ref{eq50}) holds for $p$, or equivalently, for all $0\leq j\leq d-1$ $$\lim_{n\rightarrow \infty} r(p^n)_j=\frac 1d,$$
where $r(p^n)_j$ denotes the coefficient of $x^j$ in the remainder of $p^n$ modulo $x^d-1.$
\end{corollary}
\begin{proof}Let $p$ be a polynomial satisfying the conditions of Lemma \ref{egy} and let $p(1)=1$. Let us denote the coefficient vector of $r(p)$ by $c$ and the circulant matrix $circ(c)$ by $C$. By Lemma \ref{egy} if $m\geq d-1$ then $C^m$ is a positive matrix. We use $J$ to denote the $d \times d$  matrix in which every entry is $\frac 1d$. By Lemma \ref{lem}, \hspace{0.5mm} $\displaystyle \lim_{n \rightarrow \infty}C^{(d-1)n}=J$. We may also express this fact by using the max norm of matrices: $$\lim_{n \rightarrow \infty}||C^{(d-1)n}-J||=0.$$  We show that for all $0\leq i\leq d-1$, $\displaystyle \lim_{n \rightarrow \infty}C^{(d-1)n+i}=J$. We recall that $JC^i=J$, so we have 
\begin{equation}\label{eq60}
 ||C^{(d-1)n}C^i-J||=||C^{(d-1)n}C^i-JC^i||\leq ||C^{(d-1)n}-J||\cdot||C^i||.
\end{equation} 
Because $||C^i||$ is a constant, we obtain that $\displaystyle \lim_{n \rightarrow \infty}C^{(d-1)n+i}=J$ for all $0\leq i\leq d-1$, consequently $C^m\rightarrow J$ as $m$ tends to infinity.
\end{proof}

\hspace{-4.3mm}\textbf{Application to Pascal's triangle} \hspace{1.3mm}For any $d \geq 2$ we can extract from the $n^{th}$ row of Pascal's triangle $\big(\binom n0,\ldots,\binom nn\big)$ $d$ disjoint subsequences one for each $j=0,\ldots,d-1$:
$$\binom nj, \binom {n}{j+d}, \binom {n}{j+2d},\ldots$$  It follows from Theorem \ref{e} that the $d$ different sums $\binom nj+ \binom {n}{j+d}+ \binom {n}{j+2d}+\cdots$ are asymptotically equal to $2^n/d$ as $n\rightarrow \infty$. 

\end{document}